\documentclass[conference, compsoc]{IEEEtran}
\usepackage{amscd}
\usepackage{amssymb}
\usepackage{latexsym}
\usepackage{amsfonts}
\usepackage[all]{xypic}
\usepackage[all]{xy}
\usepackage{amsmath}
\usepackage{array}
\usepackage[cp1250]{inputenc}

\usepackage{hyperref}

\ifCLASSOPTIONcompsoc
  \usepackage[nocompress]{cite}
\else
  \usepackage{cite}
\fi
\newtheorem{lem}{Lemma}[section]{\bf}{\it}
\newtheorem{prop}[lem]{Proposition}{\bf}{\it}
{\bf}{\it}
\newtheorem{thm}[lem]{Theorem}{\bf}{\it}
\newtheorem{defin}[lem]{Definition}{\bf}{\rm}
\newtheorem{rem}[lem]{Remark}{\bf}{\rm}
{\bf}{\rm}
\newtheorem{ex}[lem]{Example}{\bf}{\rm}
{\bf}{\rm}

\newtheorem{algorithm}[lem]{Algorithm}

\def\itemfakt#1{$(#1)$}

\def\Nak#1{N_{#1}}
\def\ram{R}
\def\startype{T}
\def\przysun{\vspace{-1ex}}
\def\prz{\vspace{-1ex}}
\def\przitem{\vspace{-6px}}
\def\AA{\mathbb A}

\def\DD{\mathbb D}
\def\EE{\mathbb E}

\def\MM{\mathbb M}

\def\NN{\mathbb N}

\def\RR{\mathbb R}
\def\QQ{\mathbb Q}

\def\ZZ{\mathbb Z}

\def\CO{{\mathcal{O}}}
\def\tr{{\rm tr}}
\def\sym#1{{\rm s}(#1)}
\def\Gdach#1{\check{G}_{#1}}
\def\uti{UTIM\,}
\def\utis{UTIM}
\def\infll{{\mathbf{t}}}

\def\exec{{\rm EX}}
\def\stg{\xi}
\def\wcie#1#2{\noindent
\hspace{#1pt} #2\par\vspace{0pt}}

\def\igfpos{\beta}
\def\igfposs{\gamma}

\begin{document}

\title{Two nondeterministic positive definiteness tests\\ for unidiagonal integral matrices
\thanks{Identify applicable funding agency here. If none, delete this.}}

\IEEEpubid{\begin{minipage}{\textwidth}\ \\[12pt]
 978-1-5386-0869-2/17/\$31.00 \copyright 2017 IEEE \\
 4th International Conference on Advances in Electrical Engineering 28-30 September, 2017, Dhaka, Bangladesh
\end{minipage}}

\author{\IEEEauthorblockN{Andrzej Mr\'oz\\{\bf 2016 (2019)}$^\star$}
\IEEEauthorblockA{ Faculty of Mathematics and Computer Science\\
Nicolaus Copernicus University\\
ul. Chopina 12/18, 87-100 Toru\'n, Poland\\
Email: amroz@mat.umk.pl}
}

\maketitle

\begin{abstract}
For standard algorithms verifying positive definiteness of a  matrix $A\in\mathbb{M}_n(\mathbb{R})$ based on Sylvester's criterion, the computationally pessimistic case is this when $A$ is positive definite.
We present two algorithms realizing the same task for $A\in\mathbb{M}_n(\mathbb{Z})$,  for which the case when $A$ is positive definite is the optimistic one.  The algorithms have  pessimistic computational complexities $\mathcal{O}(n^3)$ and $\mathcal{O}(n^4)$ and they rely on performing certain edge transformations, called inflations, on the edge-bipartite graph (=bigraph) $\Delta=\Delta(A)$ associated with $A$. We provide few variants of the algorithms, including Las Vegas type randomized ones with precisely described  maximal number of steps. The algorithms work very well in practice, in many cases with a better speed than the standard tests. Moreover, the algorithms yield some additional information on the properties on the quadratic form $q_A:\mathbb{Z}^n\to\mathbb{Z}$ associated with a matrix $A$.
On the other hand, our results provide an interesting example of  an application of symbolic computing methods originally developed for different purposes, with a big potential for further generalizations in matrix problems.

This is an extended version of the article \cite{Msynasc2016} in which we discussed the algorithm of the complexity $\mathcal{O}(n^4)$.
\end{abstract}

\smallskip

{\footnotesize{{$^\star$ Written in 2016 (and put in a drawer). Updated and published in arXiv in 2019. The update relies on: few remarks at the end, new test results after a slight optimisation of the codes, updates of the literature.}}}

\section{Introduction}

Fix an integer $n\geq 1$ and a matrix $A\in\MM_n(\RR)$, i.e., $A$ is a square $n\times n$ matrix with real coefficients. Recall that $A$ is {\em positive definite} if one of the following equivalent conditions holds (see \cite[Section 7.6]{Meyer}):

\begin{enumerate}
\item $q_A(x)>0$ for every non-zero $x\in\RR^n$, where $q_A:\RR^n\to\RR$
is the quadratic form associated with $A$, defined by the formula $q_A(x)=x^\tr\cdot A\cdot x$, for every (column vector) $x\in\RR^n$.
\item All eigenvalues of the matrix $\sym{A}:=\frac{1}{2}(A^\tr+A)$ are positive ($\sym{A}$ is called the {\em symmetrization} of $A$).
\item All leading principal minors of $\sym{A}$ are positive.
\item A non-positive pivot does not emerge during row reduction of $\sym{A}$ to upper triangular form with \lq\lq Type III'' operations only (i.e., operations of adding $i$th row multiplied by scalar $\alpha\in\RR$ to $j$th row, with $i<j$).
\end{enumerate}
We recall that the positive definiteness of a  matrix is an important property exploited in many branches of theoretical and applied mathematics
(see \cite{Meyer} and \cite{Minfl2016, Mgroth2016, SimSIAM2013, Draxetal95, KasjanSimson2015, BarotPena99}) and algorithmic methods related to this concept are of big importance.

Observe that the condition 4) yields quite efficient obvious algorithm for testing whether a matrix $A\in\MM_n(\RR)$ is positive definite, with pessimistic complexity $\CO(n^3)$ (as standard Gaussian elimination). Moreover, it is clear that the pessimistic case for this algorithm is the case when $A$ is positive definite (if it is not, the row reduction process can be interrupted when  first non-positive pivot emerges).

In the present paper we provide few variants of a new algorithmic positive definiteness test for so-called {\em unidiagonal} {\em triangle-integral}  matrices (certain subset of real matrices containing the set  all matrices with integer coefficients and the unity on the diagonal, see Section 2), for which the case when a matrix is positive definite is the optimistic one. The pessimistic complexities of our algorithms are $\CO(n^3)$ and $\CO(n^4)$ (see Algorithms \ref{alg:pdti} and  \ref{alg:pdtri}, and Theorem \ref{thm:main3}). However,   for some classes of matrices,  their  performance is in practice better than the standard tests as this mentioned above. The proposed algorithms can have a useful application e.g. when there is a need  to perform the test for positive definiteness for a large set of large, dense integral matrices,  for which one knows a priori that  \lq\lq most of them'' are positive definite.

Our results are an application of the inflation techniques for bigraphs from \cite{Minfl2016}, inspired by the recent studies of  Simson \cite{SimSIAM2013} (see Section 2, see also  \cite{KasjanSimson2015.3, KosFI2012, Zajac2016}). These ideas extend older known results of Ovsienko \cite{Ovsienko}, Dr\"axler et al. \cite{Draxetal95} and Barot-de la Pe\~na \cite{BarotPena99} concerning combinatorial, algebraic and algorithmic analysis of integral quadratic forms and their roots. We refer to \cite{Minfl2016} for more complete bibliography and more detailed introduction to the topic. On the other hand, inflations of bigraphs, matrices and integral quadratic forms are an important ingredient of the new lively area -- the Coxeter spectral analysis of bigraphs and their morsifications \cite{SimSIAM2013, GaSimZajLAA2016,KasjanSimson2015,KasjanSimson2015.2,KasjanSimson2015.3,MrozCE2016,
Sim2016.1,Sim2016.2,Zajac2016,Minfl2016,Mgroth2016}, cf. \cite{BarotPena99}, as well as the Coxeter type study of finite dimensional algebras  and their bounded derived categories, see
\cite{MrozPena2014, MrozPena2016, SimSIAM2013, Mgroth2016}. In the present paper we introduce slightly different than in \cite{SimSIAM2013, Minfl2016,Mgroth2016,Sim2016.1, Sim2016.2} conventions concerning an approach to interrelations between matrices and bigraphs (more elementary but at the same time more general).

The  paper can be treated as a continuation of a computer algebra approach to matrix problems presented by the author at SYNASC 2012 \cite{KasjanMroz2012}, see also \cite{DMdomestic, DMdomestic2, DMM2,DMM3,DM2, GKM}. The paper is also an extended version of \cite{Msynasc2016} where we discussed Algorithm \ref{alg:pdti} (of the complexity $\mathcal{O}(n^4)$) only.  Algorithm \ref{alg:pdtri} (of the complexity $\mathcal{O}(n^3)$)   is the main contribution of the present paper. Implementations of the algorithms discussed here are included as a part of Maple packages \cite{CIG} (see also author's homepage \cite{Mhomepage} to download different projects, including  C++ implementation of the inflation algorithm with the graphical presentation).

\section{Preliminary notions and facts}

By a {\em $($real$)$ quadratic form}  of rank $n\geq 1$ we mean a mapping $q:\RR^n\to \RR$ defined by a homogeneous polynomial of second degree
 \begin{equation}
 q(x)=\sum_{i\leq j}q_{ij}x_{i}x_{j},
 \end{equation}
 where $x=(x_1,x_2,\ldots,x_n)$ and $q_{ij}$ are fixed real numbers, for $1\leq i\leq j\leq n$. We say that $q$ is a {\em unit} quadratic form if $q_{ii}=1$, for each $i=1,\ldots,n$. If all the coefficients $q_{ij}$ are integers then
 we treat $q$ as the mapping $q:\ZZ^n\to\ZZ$  and call it an {\em integral} quadratic form. Note that  $q(x)>0$ for every non-zero $x\in\RR^n$ if and only if $q(x)>0$ for every non-zero $x\in\ZZ^n$, in case $q$ is an integral quadratic form (see \cite[Section 1.1]{Ri}).

 Each matrix $A=[a_{ij}]\in\MM_n(\RR)$
 defines a quadratic form $q_A:\RR^n\to\RR$ by the following formula
\begin{equation}\begin{array}{rcl}
 q_A(x)&:=&x^\tr\cdot A\cdot x\ \,=\ \,\sum_{i,j}a_{ij}x_{i}x_{j}\smallskip\\
 &=&\sum_{i=1}^na_{ii}x_{i}^2 + \sum_{i<j}(a_{ij}+a_{ji})x_{i}x_{j}\end{array}\label{eq:qfm}\end{equation}
 (we always treat $x\in\RR^n$ as a column vector in this context).

  Given a matrix $A\in\MM_n(\RR)$, we set $\sym{A}:=\frac{1}{2}(A^\tr+A)$; moreover, we denote by  $\nabla(A)$,   the (upper) {\em triangularisation} of  $A$, i.e., $\nabla(A)=[a^\nabla_{ij}]\in\MM_n(\RR)$ is an upper triangular matrix such that $a^\nabla_{ij}=a_{ij}+a_{ji}$ for $i<j$; $a^\nabla_{ij}=a_{ij}$ for $i=j$; and $a^\nabla_{ij}=0$ for $i>j$ (see \cite{SimSIAM2013, Mgroth2016}
 , cf. Example \ref{ex:triangbigr}).

We often use the following simple observation.

\begin{lem}\label{lem:symtriang} {\em Fix $n\geq 1$ and a matrix $A\in\MM_n(\RR)$. Then
\begin{itemize}
\item[\itemfakt{a}] $q_A=q_{\sym{A}}=q_{\nabla(A)}$;

\vskip2px

\item[\itemfakt{b}] the following conditions are equivalent:
\begin{itemize}\przitem
\item[$\bullet$] $A$ is positive definite,
\item[$\bullet$] $\sym{A}$ is positive definite,
\item[$\bullet$] $\nabla(A)$ is positive definite.
\end{itemize}
\end{itemize}
}\end{lem}

{\em Proof.} (a) Since $(x^\tr\cdot A^\tr\cdot x)^\tr=x^\tr\cdot A\cdot x\in\MM_1(\RR)=\RR$, we have
$$x^\tr\cdot \sym{A}\cdot x=\frac{1}{2}(x^\tr\cdot A^\tr\cdot x+x^\tr\cdot A\cdot x)=x^\tr\cdot A\cdot x,$$
for any $x\in\RR^n$, so $q_{\sym{A}}=q_A$. The equality $q_A=q_{\nabla(A)}$ follows from \eqref{eq:qfm}. The assertion (b) is an immediate consequence of (a) (see Introduction).\hfill$\Box$

\smallskip

We introduce the following technical notions to characterize the input matrices for our algorithms.

\begin{defin}\label{def:triang}
Fix $n\geq 1$ and a matrix $A\in\MM_n(\RR)$.
\begin{itemize}
\item[(a)] We say that $A$ is a {\em unidiagonal}\footnote{Note that the term \lq\lq unidiagonal'' appears in the literature with various different meanings.}
    matrix, if $A$ has 1's on the diagonal (equivalently,  $q_A$ is a unit quadratic form).
\item[(b)] We say that $A$ is a {\em triangle-integral} matrix, if
its triangularisation $\nabla(A)$ has integral coefficients (equivalently, $q_A$ is an integral quadratic form).
\end{itemize}
\end{defin}
Clearly, each integral matrix $A\in\MM_n(\ZZ)$ is triangle-integral.
If $A\in\MM_n(\RR)$ is a unidiagonal and triangle-integral matrix, we say that $A$ is a \utis, for short.

\smallskip

The following fact provides useful restrictions for the coefficients of a positive definite matrix $A$  in terms of its diagonal.
\begin{lem}\label{lem:posmatrix} {\em Let $A=[a_{ij}]\in\MM_n(\RR)$, $n\geq 1$, be a
positive definite matrix. Then
\begin{itemize}
\item[\itemfakt{a}] $a_{ii}>0$, for each $1\leq i\leq n$;
\item[\itemfakt{b}] $-(a_{ii}+a_{jj})<a_{ij}+a_{ji}<a_{ii}+a_{jj}$, for each $1\leq i<j\leq n$;
\item[\itemfakt{c}] if $A$ is a \uti\, then each coefficient of the triangularisation $\nabla(A)$ of $A$ belongs to the set $\{-1,0,1\}$.
\end{itemize}    }
\end{lem} 

{\em Proof}.  Since $A$ is positive definite, we have
$$\begin{array}{c}q_A(e_i)=a_{ii}>0,\smallskip\\
q_A(e_i+e_j)=a_{ii}+a_{jj}+a_{ij}+a_{ji}>0,\smallskip\\
q_A(e_i-e_j)=a_{ii}+a_{jj}-a_{ij}-a_{ji}>0,\end{array}$$
for each $i<j$, by the formula \eqref{eq:qfm}, where $e_1, \ldots, e_n$ denote the standard basis of  $\RR^n$. This implies (a) and (b).

The assertion (c) follows by (b) and Definition \ref{def:triang} (cf.  \cite[Lemma 4.1(a)]{KosFI2012}).\hfill$\Box$

\smallskip

The following specialization of a (signed) multigraph from \cite{SimSIAM2013} is an important combinatorial  tool exploited in our algorithms.

\begin{defin}\label{def:bigraph}
An {\em edge-bipartite graph} (or, a {\em bigraph}, for short),  is a pair $(\Delta_0, \Delta_1)$, where $\Delta_0$ is a finite
non-empty set of vertices  and $\Delta_1$ is a finite multiset of edges (i.e., unordered pairs of vertices) equipped with a disjoint union bipartition $\Delta_1=\Delta_1^-\cup\Delta_1^+$ such that the (multi)set $\Delta_1(a, b)=\Delta_1(b, a) \subseteq \Delta_1$
of edges connecting the vertices $a$ and $b$ splits into a disjoint
union $\Delta_1(a, b) = \Delta_1^-(a, b) \cup \Delta_1^+(a, b)$ with $\Delta_1^-(a, b) \subseteq
\Delta_1^-$  and $\Delta_1^+(a, b) \subseteq
\Delta_1^+$, for each pair of
vertices $a,b\in\Delta_0$, and either $\Delta_1(a, b) = \Delta_1^-(a,b)$ or $\Delta_1(a, b) = \Delta_1^+(a,b)$.
 \end{defin}

In other words, a bigraph is an undirected multigraph  equipped with two kinds of edges: these from
$\Delta_1^-$ (called  {\em solid edges}) and  $\Delta_1^+$ (called {\em dotted edges}); and two vertices can be connected by only one kind of these two kinds of edges. We define a bigraph $\Delta=(\Delta_0,\Delta_1)$ to be loop-free if $\Delta$ has no loops, that is, $\Delta_1(a,a)=\emptyset$ for each vertex $a\in\Delta_0$. We say that a bigraph $\Delta$ is connected, if the underlying graph $\overline{\Delta}$, obtained from $\Delta$ by replacing all
dotted edges  by the continuous ones  is connected. Note that a bigraph $\Delta=(\Delta_0,\Delta_1)$ with $\Delta_1^+=\emptyset$, that is, a bigraph  without dotted edges, is a usual (multi)graph; and each undirected (multi)graph $\Delta$ can be viewed as a bigraph in this way. Later on, we always assume that each  bigraph $\Delta=(\Delta_0,\Delta_1)$ with $|\Delta_0|=n\geq 1$  has fixed labeling of  vertices by first $n$ natural numbers, i.e., $\Delta_0=\{1,2,\ldots,n\}$.

Given $n\geq 1$ and a \uti $A\in\MM_n(\RR)$, we associate with $A$ the (unique) loop-free bigraph $\Delta:=\Delta(A)=(\Delta_0,\Delta_1)$ such that  \begin{itemize}
\item $\Delta_0=\{1,2,\ldots,n\}$;
\item $|\Delta_1^-(i,j)|=|a_{ij}^\nabla|$ if $a_{ij}^\nabla<0$ and $i<j$;
\item $|\Delta_1^+(i,j)|=a_{ij}^\nabla$ if $a_{ij}^\nabla>0$ and $i<j$,
\end{itemize}
where $\nabla(A)=[a^\nabla_{ij}]\in\MM_n(\RR)$ is the {triangularisation} of  $A$. We say that a \uti $A\in\MM_n(\RR)$ is {\em connected} if the bigraph $\Delta(A)$ is connected.

\begin{rem}\label{rem:loops}
(a) In case a triangle-integral matrix $A\in\MM_n(\RR)$ is not unidiagonal one can associate with $A$ certain bigraph with loops (cf. \cite{Mgroth2016,KasjanSimson2015,KasjanSimson2015.2, KasjanSimson2015.3,Sim2016.2}). However,  in this paper we restrict to unidiagonal matrices and loop-free bigraphs. The general case is technically  more complex and more delicate; we plan to study it in the future paper.

(b) Lemma \ref{lem:posmatrix}(c) implies that if a \uti $A$ is positive definite then the bigraph $\Delta(A)$ does not have multiple edges.

(c) The notion of a triangle-integral matrix we introduced above can be treated as a slight generalization of Simson's concept of a {\em rational morsification} (see \cite[Section 3]{Sim2016.1}, cf. \cite{SimSIAM2013, Mgroth2016}). More precisely, one checks that if $A\in\MM_n(\QQ)$ is a rational morsification of a bigraph $\Delta$, then $A$ is a triangle-integral matrix. Moreover, if $\Delta$ is loop-free, then $A$ is unidiagonal (see \cite{Sim2016.1} for the details).
\end{rem}

Note that $\Delta(A)=\Delta(A')$ iff $\nabla(A)=\nabla(A')$, for a pair of \utis's $A,A'\in\MM_n(\RR)$. Moreover, given a
loop-free bigraph $\Delta=(\Delta_0,\Delta_1)$ with $\Delta_0=\{1,\ldots,n\}$, there exists a unique upper-triangular unidiagonal matrix $\Gdach{\Delta}\in\MM_n(\ZZ)$ such that $\Delta(\Gdach{\Delta})=\Delta$. We call $\Gdach{\Delta}$ the (upper-triangular) {\em Gram matrix} of $\Delta$ (this definition is equivalent to Simson's original definition of the upper-triangular Gram matrix of $\Delta$  from \cite{SimSIAM2013}). Note that if $\Delta=\Delta(A)$, for a \uti $A\in\MM_n(\RR)$, then $\Gdach{\Delta}=\nabla(A)$. We say that a loop-free bigraph $\Delta$ is {\em positive} (and we write $\Delta>0$) if the Gram matrix $\Gdach{\Delta}$ is positive definite.

\begin{ex}\label{ex:triangbigr}
Let $A\in\MM_4(\RR)$ be the following \uti
$$A=\left[ \begin {array}{cccc} 1&-\frac{1}{2}&1+\sqrt {2}&1\\ -\frac{3}{2}&1&0&0\\ -\sqrt {2}&1&1&7\\ -2
&0&-5&1\end {array} \right].$$
Then its triangularization $\nabla(A)$ and the associated bigraph $\Delta=\Delta(A)$ look as follows:
$$\nabla(A)=\left[ \begin {array}{rrrr} 1&-2&1&-1\\ 0&1&1&0
\\ 0&0&1&2\\ 0&0&0&1\end {array}
 \right]\ \ \Delta:{\scriptsize
\xymatrixrowsep{12pt}
\xymatrixcolsep{12pt}
{\xymatrix{
 &1\ar@/^3px/@{-}[rd]\ar@/_3px/@{-}[rd]\ar@{.}[dd]\ar@{-}[ld]& \\
4\ar@/^3px/@{.}[rd]\ar@/_3px/@{.}[rd]& &2\ar@{.}[ld]\\
 &3& \\}}}$$
 Note that $\Delta$ is not positive (equivalently, $A$ is not positive definite, see Lemma \ref{lem:symtriang}), since $\Delta$ has multiple edges (see Remark \ref{rem:loops}(b)).
\end{ex}

\section{Generalities on inflations of bigraphs}

We recall from  \cite{SimSIAM2013, KosFI2012} and \cite{Minfl2016} the following concepts of inflations of bigraphs.

 \begin{defin}\label{def:infl} Fix an integer $n\geq 1$ and a loop-free bigraph $\Delta=(\Delta_0,\Delta_1)$ with $\Delta_0=\{1,\ldots,n\}$ and the Gram matrix $$\Gdach{\Delta}={\scriptsize\left[\begin{array}{ccccc}
     1&d_{1,2}^\Delta&\cdots&d_{1,n\!-\!1}^\Delta&d_{1,n}^\Delta\\
     0&1&\cdots&d_{2,n\!-\!1}^\Delta&d_{2,n}^\Delta\\
     \vdots&\vdots& \ddots&&\vdots\\
     0&0&\cdots&1&d_{n\!-\!1,n}^\Delta\\
     0&0&\cdots&0&1
     \end{array}\right]}\in\MM_n(\ZZ).$$
 \begin{itemize}
\item[{\rm (a)}] An {\em inflation at a vertex} $a\in\Delta_0$ is an operation associating with $\Delta$ the new bigraph $\hat{\Delta}:=\infll_a\Delta$ such that $\hat{\Delta}_0:=\Delta_0=\{1,\ldots,n\}$ and the edges of $\hat{\Delta}$ are obtained from $\Delta$ by replacing any solid edge (resp. dotted edge) in $\Delta_1$ incident with $a\in\Delta_0$ by a dotted one (resp. solid one) in $\hat{\Delta}_1$; the remaining edges in $\Delta$ stay unchanged in $\hat{\Delta}$.
\item[{\rm (b)}] Given a pair $a,b\in\Delta_0$ of distinct vertices such that
the set $\Delta_1^+(a,b)$ of dotted edges between $a$ and $b$ is non-empty, an {\em inflation at the pair} $(a,b)$ is an operation associating with $\Delta$ the new bigraph $\hat{\Delta}:=\infll_{a,b}\Delta$ with $\hat{\Delta}_0:=\Delta_0=\{1,\ldots,n\}$ and the edges of $\hat{\Delta}$  defined as follows:
    \begin{itemize}\przitem
    \item[$\bullet$] replace each of the dotted edges in $\Delta_1^+(a,b)$ by a solid one  in $\hat{\Delta}_1^-(a,b)$;

    \item[$\bullet$] the multiset of edges $\hat{\Delta}_1(b,c)=\hat{\Delta}_1(c,b)$ is defined by setting $d^{\hat{\Delta}}_{b,c} := d^\Delta_{bc}-d^\Delta_{ac}d^\Delta_{ab}$,\, for each $c\neq a,b$, where we set $d^\Delta_{ij}:=d^\Delta_{ji}$ and $d^{\hat{\Delta}}_{ij}:=d^{\hat{\Delta}}_{ji}$ if $i>j$;

    \item[$\bullet$] each of the remaining edges in $\Delta$ stay unchanged in $\hat{\Delta}$.
    \end{itemize}
\end{itemize}
\end{defin}
Note that $\infll_{a,b}\Delta$ is defined if and only if so is $\infll_{b,a}\Delta$, but $\infll_{a,b}\Delta\neq \infll_{b,a}\Delta$ in general (see \cite[Example 4.5(a)]{Minfl2016}). Observe that  a bigraph $\Delta$ is a (multi)graph (i.e., $\Delta$ does not have dotted edges) if and only if none of inflations at a pair is defined for $\Delta$.

\begin{ex}\label{ex:infl}
Let $\Delta$ be the following bigraph:
$${\scriptsize
\xymatrixrowsep{6pt}
\xymatrixcolsep{4pt}
\Delta:
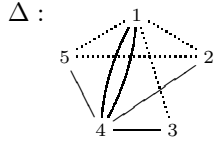
{\xymatrix{
 & &1\ar@{.}[lld]\ar@/_3px/@{-}[lddd]\ar@/^3px/@{-}[lddd]\ar@{.}[rddd]\ar@{.}[rrd]& & \\
5\ar@{.}[rrrr]\ar@{-}[rdd]& & & &2\\
 & & & & \\
 &4\ar@{-}[rrruu]\ar@{-}[rr]& &3& \\
}}}$$

\noindent
Then $\infll_{2,1}\Delta$ and $\infll_{1}\Delta$ look as follows
$$\infll_{2,1}\Delta:{\scriptsize
\xymatrixrowsep{6pt}
\xymatrixcolsep{4pt}
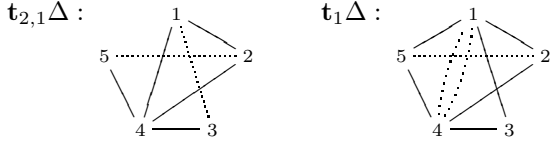
{\xymatrix{
 & &1\ar@{-}[lddd]\ar@{.}[rddd]\ar@{-}[rrd]& & \\
5\ar@{.}[rrrr]\ar@{-}[rdd]& & & &2\\
 & & & & \\
 &4\ar@{-}[rrruu]\ar@{-}[rr]& &3& \\
}}}\hspace{0.8cm} \infll_{1}\Delta:{\scriptsize
\xymatrixrowsep{6pt}
\xymatrixcolsep{4pt}
{\xymatrix{
 & &1\ar@{-}[lld]\ar@/_3px/@{.}[lddd]\ar@/^3px/@{.}[lddd]\ar@{-}[rddd]\ar@{-}[rrd]& & \\
5\ar@{.}[rrrr]\ar@{-}[rdd]& & & &2\\
 & & & & \\
 &4\ar@{-}[rrruu]\ar@{-}[rr]& &3& \\
}}}$$
\end{ex}
We refer to \cite{Minfl2016} and \cite{KosFI2012, SimSIAM2013} for more examples and the discussion on other properties of inflations.

Following \cite{Minfl2016} we formalize an execution of inflation algorithm as follows.
\begin{defin}\label{def:exec} Fix a loop-free bigraph $\Delta$.
 The sequence $\exec:=(\infll^{(s)},\infll^{(s-1)},\ldots, \infll^{(1)}$) is called an {\em execution of inflation algorithm on} $\Delta$ {\em resulting in} $\hat{\Delta}$ if
\begin{itemize}
  \item
    each of its terms $\infll^{(j)}$ has the form $\infll^{(j)}=\infll_c$, with $c\in\Delta_0$, or $\infll^{(j)}=\infll_{a,b}$, with $a\neq b\in\Delta_0$,
  \item $\infll^{(j)}\infll^{(j-1)}\cdots \infll^{(1)}\Delta$ is defined, for any $j=1,\ldots,s$,
   \item $\infll^{(s)}\infll^{(s-1)}\cdots \infll^{(1)}\Delta=\hat{\Delta}$.
\end{itemize}

The integer $s\geq 1$ is called the {\em length} of the execution $\exec$.
We admit also an empty execution $\exec=\emptyset$ on $\Delta$ of length
$s=0$, resulting in $\hat{\Delta}=\Delta$.
 \end{defin}

Clearly, the inflation algorithm in the sense of Definition \ref{def:exec} is not deterministic (see also Section 4 and \cite{Minfl2016}).

\begin{prop}\label{prop:preserve} {\em Let $\Delta$ be a loop-free bigraph and $\exec$ an  execution of inflation algorithm on $\Delta$ resulting in $\hat{\Delta}$.
\begin{itemize}
\item[\itemfakt{a}] If $\Delta$ is connected then so is $\hat{\Delta}$;
\item[\itemfakt{b}]  $\Delta$ is  positive if and only if so is $\hat{\Delta}$.
\end{itemize}
}
\end{prop}

{\em Proof.} Both  assertions follow from \cite[Corollary 3.4]{Minfl2016} (see also \cite[Lemma 3.3]{Minfl2016}, cf. \cite[Lemma 4.1]{KosFI2012}).\hfill$\Box$

\smallskip

It was proved in \cite[Theorem 3.1]{SimSIAM2013} (see also \cite[Theorem 4.4]{Minfl2016} and \cite{KosFI2012}) that a connected loop-free bigraph $\Delta$ is positive if and only if there exists an execution $\exec$ on $\Delta$ resulting in one of the following {\em Dynkin graphs} $D$ (see \cite{SimSIAM2013, Minfl2016}):
\begin{center}
\begin{tabular}{cc}
$\AA_n$&
\begin{picture}(135,10)
\put(-1,-1){{\scriptsize $\bullet$}} \put(20,-1){{\scriptsize $\bullet$}} \put(44,-1){$\cdots$}
\put(80,-1){{\scriptsize $\bullet$}} \put(100,-1){{\scriptsize $\bullet$}}
\put(0,1){\line(1,0){40}} \put(60,1){\line(1,0){40}}
\end{picture}\\
 & \hspace{-1cm}($n$ vertices, $n\geq 1$)  \\
$\DD_n$&
\begin{picture}(135,30)
\put(-1,-1){{\scriptsize $\bullet$}} \put(20,-1){{\scriptsize $\bullet$}} \put(40,-1){{\scriptsize
$\bullet$}} \put(22,1){\line(0,1){20}} \put(20,20){\scriptsize $\bullet$} \put(64,-1){$\cdots$}
\put(100,-1){{\scriptsize $\bullet$}} \put(120,-1){{\scriptsize $\bullet$}}
\put(0,1){\line(1,0){60}} \put(80,1){\line(1,0){40}}
\end{picture}\\
 & \hspace{-1cm}($n$ vertices, $n\geq 4$)  \\
$\EE_6$&
\begin{picture}(135,30)
\put(-1,-1){{\scriptsize $\bullet$}} \put(20,-1){{\scriptsize $\bullet$}} \put(40,-1){{\scriptsize
$\bullet$}} \put(42,1){\line(0,1){20}} \put(40,20){\scriptsize $\bullet$} \put(60,-1){{\scriptsize
$\bullet$}} \put(80,-1){{\scriptsize $\bullet$}} \put(0,1){\line(1,0){80}}
\end{picture}\\
$\EE_7$&
\begin{picture}(135,30)
\put(-1,-1){{\scriptsize $\bullet$}} \put(20,-1){{\scriptsize $\bullet$}} \put(40,-1){{\scriptsize
$\bullet$}} \put(42,1){\line(0,1){20}} \put(40,20){\scriptsize $\bullet$} \put(60,-1){{\scriptsize
$\bullet$}} \put(80,-1){{\scriptsize $\bullet$}} \put(100,-1){{\scriptsize $\bullet$}}
\put(0,1){\line(1,0){100}}
\end{picture}\\
$\EE_8$&
\begin{picture}(135,30)
\put(-1,-1){{\scriptsize $\bullet$}} \put(20,-1){{\scriptsize $\bullet$}} \put(40,-1){{\scriptsize
$\bullet$}} \put(42,1){\line(0,1){20}} \put(40,20){\scriptsize $\bullet$} \put(60,-1){{\scriptsize
$\bullet$}} \put(80,-1){{\scriptsize $\bullet$}} \put(100,-1){{\scriptsize $\bullet$}}
\put(120,-1){{\scriptsize $\bullet$}} \put(0,1){\line(1,0){120}}
\end{picture}\\

\end{tabular}

\end{center}
In \cite[Theorem 4.4]{Minfl2016} we enhanced this result by describing precise bounds for the lengths of few variants of such executions on positive bigraphs. The further discussion of the paper, concerning the announced positive definiteness tests, is based on some conclusions from \cite{Minfl2016} (cf. \cite[Remark 4.6(b)]{Minfl2016}).

\section{Executions of inflation algorithm}

One of the crucial steps in our algorithmic solutions is
the following \lq\lq generic'' inflation procedure, which is a specialization of randomized Algorithm 4.3 from \cite{Minfl2016}:

\begin{algorithm}\label{alg:random}

\smallskip

\noindent
 {\bf Input:} a loop-free bigraph $\Delta$ with $n\geq 1$ vertices and an integer $\stg\in\{0,1,2,3\}$ (called a {\em strategy}).

\smallskip

\noindent
{\bf Output (if it stops):} a loop-free (multi)graph $D$  with $n\geq 1$ vertices such that $D$ is connected if so is $\Delta$, and $D>0$ if and only if $\Delta>0$.

\medskip

\noindent {\tt InflationsAtPair}$(\Delta,\, \stg)$

\smallskip

\noindent
\begin{tabular}{rl}
{\tt 1}&\wcie{0}{{\tt while}\,\, $\Delta$\ \ {\tt has dotted edge do \{}}\\
{\tt 2}&\wcie{10}{{\tt select dotted edge}\, $a\cdots b$}\\
{\tt 3}&\wcie{20}{{\tt with respect to strategy} $\stg${\tt ;}}\\
{\tt 4}&\wcie{10}{{\tt set} $\Delta := \,\infll_{a,b}\Delta${\tt ;}}\\
{\tt 5}&\wcie{0}{{\tt \}}}\\
{\tt 6}&\wcie{0}{{\tt return}\, $D := \Delta${\tt ;}}\\
\end{tabular}

\end{algorithm}

The strategy $\stg$ can be one of the following:
\begin{center}
\begin{tabular}{c|l}
$\stg$&strategy for lines {\tt 2-3} of Algorithm \ref{alg:random}\\
\hline
0 & find first dotted edge\\
1 & find last dotted edge\\
2 & find first or last dotted edge, randomly\\
3 & find a random dotted edge\\\hline
\end{tabular}
\end{center}
By first (resp. last) dotted edge above we mean a
dotted edge $a\cdots b$ such that the pair $(a,b)$ is the smallest (resp. the greatest) element in $\NN^2$ with respect to the lexicographic order, among all dotted edges.

In general Algorithm \ref{alg:random} does not have a stop property, but its partial correctness, i.e., the properties of the output graph $D$ in case the algorithm stops, follow from Proposition \ref{prop:preserve}. Note that in case the strategy $\stg$ equals $0$ or $1$ (resp. $2$ or $3$) then Algorithm \ref{alg:random} is a deterministic (resp. nondeterministic randomized) algorithm.

Fix an integer $n\geq 1$ and an integral quadratic form $q:\ZZ^n\to\ZZ$. A vector  $v\in\ZZ^n$ is called  a {\em root} of $q$ if $q(v)=1$ (see \cite{Minfl2016, Ri}). We say that a vector $v=[v_1,\ldots,v_n]^\tr\in\ZZ^n$ is {\em positive} (resp. {\em sincere}) if $v\neq 0$ and $v_i\geq 0$ (resp. $v_i\neq 0$), for all $i=1,\ldots,n$.

The following conclusion from \cite[Theorem 4.4]{Minfl2016} provides the most important properties of Algorithm \ref{alg:random} in context of our goals.

\begin{thm}\label{thm:main1}
{\em Let $\Delta$ be a connected loop-free bigraph with $n\geq 1$ vertices and  $\stg\in\{0,1,2,3\}$, an  arbitrary fixed strategy. Then the following hold:
\begin{itemize}
\item[\itemfakt{a}] The bigraph $\Delta$ is positive if and only if
 Algorithm \ref{alg:random}, applied for $\Delta$ and $\stg$, stops after performing at most $\igfpos(n)$ inflations, and the returned graph $D$ is a Dynkin graph, where
 $$\igfpos(n):={\scriptsize\left\{\begin{array}{cl}
        \frac{1}{2}(n^2-n),& n\in\{1,2,3\},\\
        30,& n=6,\\
        56,& n=7,\\
        112,& n=8,\\
        n^2-2n,& \mbox{elsewhere}.
        \end{array}\right.}$$
\item[\itemfakt{b}] Assume that the integral quadratic form $q_\Delta:=q_{\Gdach{\Delta}}:\ZZ^n\to\ZZ$ admits a positive sincere root. Then  $\Delta$ is positive if and only if
 Algorithm \ref{alg:random}, applied for $\Delta$ and $\stg$,  stops after performing at most $\igfposs(n)$ inflations, and the returned graph $D$ is a Dynkin graph, where
 $$\igfposs(n):={\scriptsize\left\{\begin{array}{cl}
        0,& n\in\{1,2,3\},\\
        5,& n=6,\\
        10,& n=7,\\
        21,& n=8,\\
        n-3,& \mbox{elsewhere}.
        \end{array}\right.}$$
\end{itemize}
}
\end{thm}

{\em Proof.} The assertion (a) (resp. (b)) follows from \cite[Theorem 4.4(d)]{Minfl2016} (resp. \cite[Theorem 4.4(b)]{Minfl2016}). We only recall that the proof of \cite[Theorem 4.4]{Minfl2016} relies on a careful analysis of changes of  root systems of bigraphs under inflations, and known properties of roots of Dynkin graphs (cf. \cite{Ovsienko, SimSIAM2013, KosFI2012}).\hfill$\Box$

\smallskip

Observe that Theorem \ref{thm:main1} provides bounds for the number of inflations depending only on the number of vertices $n$ of $\Delta$, independently on the chosen strategy $\stg$. In particular, it shows that  Algorithm \ref{alg:random} is a Las Vegas type randomized algorithm in case $\stg\in\{2,3\}$ and $\Delta>0$.

\begin{algorithm}\label{alg:psr}

\smallskip

\noindent
 {\bf Input:} a loop-free connected bigraph $\Delta=(\Delta_0,\Delta_1)$ with $\Delta_0=\{1,\ldots,n\}$, $n\geq 1$.

\smallskip

\noindent
{\bf Output:} a loop-free connected bigraph $\hat{\Delta}$  with $n\geq 1$ vertices such  the quadratic form $q_{\hat{\Delta}}:\ZZ^n\to\ZZ$ admits a positive sincere root, and $\Delta>0$ if and only $\hat{\Delta}>0$.

\medskip

\noindent {\tt InflationsToPosSincereRoot}$(\Delta)$

\smallskip

\noindent
\begin{tabular}{rl}
{\tt 1}&\wcie{0}{{\tt set} $S:=\{1\}${\tt ;}}\\
{\tt 2}&\wcie{0}{{\tt while}\,\, $|S|<n$\ \ {\tt do \{}}\\
{\tt 3}&\wcie{10}{{\tt choose $(a, b)\in S\times(\Delta_0\setminus S)$}}\\
{\tt 4}&\wcie{20}{{\tt such that $\Delta_1(a,b)\neq \emptyset$;}}\\
{\tt 5}&\wcie{10}{{\tt if $\Delta_1^+(a,b)=\emptyset$ then set $\Delta:=\infll_b\Delta$;}}\\
{\tt 6}&\wcie{10}{{\tt  set $\Delta:=\infll_{b,a}\Delta$;}}\\
{\tt 7}&\wcie{10}{{\tt  set $S:=S \cup\{b\}$;}}\\
{\tt 8}&\wcie{0}{{\tt \}}}\\
{\tt 9}&\wcie{0}{{\tt return}\, $\hat{\Delta}:=\Delta${\tt ;}}\\
\end{tabular}
\end{algorithm}

\begin{thm}\label{thm:main2}
{\em Algorithm \ref{alg:psr} is correct and applied to a loop-free connected bigraph $\Delta$ with $n\geq 1$ vertices, it performs at most $n-1$ inflations at a vertex and at most $n-1$ inflations at a pair.}
\end{thm}

{\em Proof.} First note that a pair $(a, b)\in S\times(\Delta_0\setminus S)$ such that $\Delta_1(a,b)\neq \emptyset$ as in lines {\tt 3-4} always exists by the connectedness of $\Delta$ (see Proposition \ref{prop:preserve}(a)). Next,
the instruction in the line {\tt 5} guarantees that the inflation $\infll_{b,a}$ in  {\tt 6} is defined. The stop property follows from the line {\tt 7} and the construction of the loop.

Let $\hat{\Delta}$ be the bigraph returned by the algorithm applied for $\Delta$.
Proposition \ref{prop:preserve}(b) implies that $\Delta>0$ if and only if $\hat{\Delta}>0$. For the proof of the fact that $q_{\hat{\Delta}}$ admits a positive sincere root $v\in\ZZ^n$ we refer to \cite[Lemma 3.3]{Minfl2016} and the proof of \cite[Theorem 4.4(c)]{Minfl2016} (we only outline the main idea:  $v$ is constructed from the \lq\lq trivial'' root $e_1=[1,0,\ldots,0]^\tr$ of the input bigraph $\Delta$ by consecutive changes induced by  inflations from the line {\tt 6}).

The numbers of performed inflations at a vertex and at a pair follow obviously from the construction of the main loop (note that $|\Delta_0\setminus\{1\}|=n-1$).\hfill$\Box$

\section{Positive definiteness tests by inflations}

We use the notation from the previous section. Let {\tt InflationsAtPairPos} (resp. {\tt InflationsAtPairPoss}) denotes the procedure which works exactly as
{\tt InflationsAtPair} for a bigraph $\Delta$ with $n\geq 1$ vertices and a strategy $\stg$, but breaks the main loop after performing $\igfpos(n)$ (resp. $\igfposs(n)$) inflations from the line {\tt 4}, see Theorem \ref{thm:main1}.

Now we have everything to formulate the main algorithms of the paper. Note that for the simplicity of the presentation we restrict to connected \utis's. Extending these ideas  for the general case (i.e., for arbitrary \utis) is not a very hard task.

\begin{algorithm}\label{alg:pdti}

\smallskip

\noindent
 {\bf Input:} a connected \uti $A\in\MM_n(\RR)$ with  $n\geq 1$, and a strategy $\stg\in\{0,1,2,3\}$ (as in Algorithm \ref{alg:random}).

\smallskip

\noindent
{\bf Output:} {\tt true} if $A$ is positive definite, or {\tt false} otherwise.

\medskip

\noindent {\tt PosDefTestByInflations}$(A,\, \stg)$

\smallskip

\noindent
\begin{tabular}{rl}
{\tt 1}&\wcie{0}{{\tt set}\, $\Delta := \Delta(A)${\tt;}}\\
{\tt 2}&\wcie{0}{{\tt set}\, $D := $ {\tt InflationsAtPairPos}$(\Delta,\, \stg)${\tt ;}}\\
{\tt 3}&\wcie{0}{{\tt if}\, $D$ is a Dynkin graph\, {\tt then} {\tt return}\, {\tt true}{\tt ;}}\\
{\tt 4}&\wcie{0}{{\tt else}\,  {\tt return}\, {\tt false}{\tt ;}}\\
\end{tabular}

\end{algorithm}

Correctness of the algorithm (for arbitrary strategy $\stg$) follows from  Theorem \ref{thm:main1}(a) and the properties of Algorithm \ref{alg:random} (see also Lemma \ref{lem:symtriang}).

\begin{algorithm}\label{alg:pdtri}

\smallskip

\noindent
 {\bf Input:} a connected \uti $A\in\MM_n(\RR)$ with  $n\geq 1$,  and a strategy $\stg\in\{0,1,2,3\}$ (as in Algorithm \ref{alg:random}).

\smallskip

\noindent
{\bf Output:} {\tt true} if $A$ is positive definite, or {\tt false} otherwise.

\medskip

\noindent {\tt PosDefTestByRootInflations}$(A,\, \stg)$

\smallskip

\noindent
\begin{tabular}{rl}
{\tt 1}&\wcie{0}{{\tt set}\, $\Delta := \Delta(A)${\tt;}}\\
{\tt 2}&\wcie{0}{{\tt set}\, $\hat{\Delta} := $ {\tt InflationsToPosSincereRoot}$(\Delta)${\tt ;}}\\
{\tt 3}&\wcie{0}{{\tt set}\, $D := $ {\tt InflationsAtPairPoss}$(\hat{\Delta},\, \stg)${\tt ;}}\\
{\tt 4}&\wcie{0}{{\tt if}\, $D$ is a Dynkin graph\, {\tt then} {\tt return}\, {\tt true}{\tt ;}}\\
{\tt 5}&\wcie{0}{{\tt else}\,  {\tt return}\, {\tt false}{\tt ;}}\\
\end{tabular}

\end{algorithm}

Correctness of the algorithm (for arbitrary strategy $\stg$) follows from  Theorem \ref{thm:main1}(b) and the properties of Algorithm \ref{alg:random} and Algorithm \ref{alg:psr} (see also Lemma \ref{lem:symtriang}).

\begin{rem}\label{rem:implem} (a) A simple optimization of Algorithm \ref{alg:pdtri} is possible. Namely, a graph $D$ (i.e., a bigraph without dotted edges) may appear already during the execution
of {\tt InflationsToPosSincereRoot}$(\Delta)$ in the line {\tt 2}. If this is the case, the algorithm can break the execution of the line {\tt 2} and jump (with $D$) to the line {\tt 4}.

\smallskip

(b) Depending on applications and implementation purposes, one can add at the beginning of Algorithms \ref{alg:pdti} and \ref{alg:pdtri}
  the test if $A$ satisfies the trivial necessary condition for positive definiteness from Lemma \ref{lem:posmatrix}(c), i.e.,
if the coefficients of the triangularization $\nabla(A)$ of $A$ belong to the set $\{-1,0,1\}$.

\smallskip

(c) From the construction of the loop in {\tt InflationsAtPair}  it is clear, that the optimistic case for Algorithms \ref{alg:pdti} and \ref{alg:pdtri} is the case when a graph $D$ emerges much before the maximal number of inflations $\igfpos(n)$ (or $\igfposs(n)$) is achieved. And this happens most often when the input matrix $A$ is positive definite (see Theorem \ref{thm:main1}).

\smallskip
(d)
Implementations of all the procedures presented in the paper are available in Maple packages \cite{CIG} (under the same names as the pseudocodes).

\smallskip
(e) An important side effect of Algorithms \ref{alg:pdti} and \ref{alg:pdtri} is the computation of the Dynkin graph $D$ for a positive definite \uti $A$. The graph $D$ is the {\em Dynkin type} of the integral quadratic form $q_A$ associated with $A$ (cf. \cite{BarotPena99, Sim2016.1}). It provides certain non-trivial additional information on $q_A$, e.g., on its root systems (see \cite{BarotPena99, Sim2016.1, SimSIAM2013, Minfl2016, Mgroth2016} for more details).
\end{rem}

In the rest of the section we prove that the pessimistic complexity of Algorithm \ref{alg:pdti} (resp. Algorithm \ref{alg:pdtri}) is $\CO(n^4)$ (resp. $\CO(n^3)$), independently on the strategy $\stg$. In particular, this  shows that Algorithms \ref{alg:pdti} and  \ref{alg:pdtri} are  Las Vegas type randomized algorithms in case $\stg\in\{2,3\}$, for arbitrary connected \uti $A$.

We start with the following technical fact, which also provides hints for implementations of the algorithms. Given a bigraph $\Delta$ with $n\geq 1$ vertices, we operate on its non-symmetric Gram matrix $\Gdach{\Delta}\in\MM_n(\ZZ)$; we count the complexity with respect to arithmetic operations performed on the coefficients of $\Gdach{\Delta}$.

\begin{prop}\label{prop:cpx}
{\em Let $\Delta=(\Delta_0,\Delta_1)$ be a loop-free bigraph with $\Delta_0=\{1,\ldots,n\}$, $n\geq 1$, and the Gram matrix $\Gdach{\Delta}\in\MM_n(\ZZ)$.
\begin{itemize}
\item[\itemfakt{a}] Inflation $\infll_a\Delta$, for $a\in\Delta_0$, can be performed with $\CO(n)$ operations.
\item[\itemfakt{b}] Inflation $\infll_{a,b}\Delta$, for $a,b\in\Delta_0$, can be performed with $\CO(n)$ operations.
\item[\itemfakt{c}] Test whether $\Delta$ is a Dynkin graph can be performed with $\CO(n^2)$ operations.
\end{itemize}
}
\end{prop}

{\em Proof.} The assertions (a) and (b) follow easily from Definition \ref{def:infl} (in both cases we change coefficients of $\Gdach{\Delta}$ corresponding to edges incident with one fixed vertex).

(c) The test whether $\Delta$ is a Dynkin graph can be performed by applying the following steps:
\begin{enumerate}
\item check if $\Delta$ is a simple graph, i.e., if $d_{ij}^\Delta\in\{0,-1\}$, for every $1\leq i<j\leq n$, where $d_{ij}^\Delta$ are the coefficients of $\Gdach{\Delta}$ as in Definition \ref{def:infl};
\item check if $\Delta$ has precisely $n-1$ edges;
\item check if $\Delta$ is connected (e.g., by applying Depth First Search);
\item find the list $\ram$ of ramifications (= vertices of degree $\geq 3$);
\item if $\ram=\emptyset$ then $\Delta=\AA_n$; if ($|\ram|>1$) or ($\ram=\{s\}$ and ${\rm deg}(s)>3$) then $\Delta$ is not a Dynkin graph;
\item assume that $\ram=\{s\}$ and ${\rm deg}(s)=3$ (i.e., $\Delta$ is a star graph with three arms); find the star type $\startype$ of $\Delta$ (i.e., $\startype$ is a list of lengths of all three paths outgoing from $s$):
    \begin{enumerate}
    \item if $\startype=(1,1,n-3)$ then $\Delta=\DD_n$;
    \item if $\startype=(1,2,2)$ then $\Delta=\EE_6$;
    \item if $\startype=(1,2,3)$ then $\Delta=\EE_7$;
    \item if $\startype=(1,2,4)$ then $\Delta=\EE_8$;
    \item otherwise, $\Delta$ is not a Dynkin graph.
    \end{enumerate}
\end{enumerate}
It is easy to observe that each of these 6 steps can be realized with the complexity at most $\CO(n^2)$.\hfill$\Box$

\smallskip

Note that the assertion (c) of the proposition is a special case of the (bi)graph isomorphism problem, for which the polynomial-time general algorithm is not known. Moreover, carefully implemented Dynkin graph recognition test based on the above 6 steps is very efficient, it appears that its performance time is negligible in Algorithms \ref{alg:pdti} and \ref{alg:pdtri}.

\begin{thm}\label{thm:main3} {\em Let $A\in\MM_n(\RR)$, $n\geq 1$, be a connected \uti and $\stg\in\{0,1,2,3\}$, an arbitrary fixed strategy. Then the following hold:
\begin{itemize}
\item[\itemfakt{a}] The pessimistic complexity of Algorithm \ref{alg:pdti}, applied to $A$ and $\stg$, is $\CO(n^4)$.
\item[\itemfakt{b}] The pessimistic complexity of  Algorithm \ref{alg:pdtri}, applied to $A$ and $\stg$, is $\CO(n^3)$.
\end{itemize}    }
\end{thm}

{\em Proof.} (a) Each execution of the body of the loop in {\tt InflationsAtPair} relies on: finding a dotted edge with respect to the strategy $\stg$, which costs $\CO(n^2)$ operations;  applying the inflation $\infll_{a,b}$, which costs $\CO(n)$ by Proposition \ref{prop:cpx}(b). The loop body is executed at most $\CO(\igfpos(n))=\CO(n^2)$ times (see Theorem \ref{thm:main1}(a)). Therefore, the execution of the line {\tt 2} in Algorithm \ref{alg:pdti} costs $\CO(n^2(n^2+n))=\CO(n^4)$ operations. Dynkin graph recognition test in the line {\tt 3} costs $\CO(n^2)$ by Proposition \ref{prop:cpx}(c), hence it does not affect the general complexity.

(b) First we estimate the complexity of the line {\tt 2} of Algorithm \ref{alg:pdtri}, i.e., the execution of {\tt InflationsToPosSincereRoot}.
By Theorem \ref{thm:main2} it performs at most $2n-2$ inflations, and this costs $\CO(n^2)$ operations by Proposition \ref{prop:cpx}. Additionally, Algorithm \ref{alg:psr} performs $n-1$ times the  edge search in its lines {\tt 3-4}, each of them of the cost $\CO(n^2)$. Therefore, the complexity of the line {\tt 2} of Algorithm \ref{alg:pdtri} is $\CO(n^2 + (n-1)n^2)=\CO(n^3)$.

The execution of {\tt InflationsAtPairPoss} in the line {\tt 3} costs $\CO(\igfposs(n)\cdot n^2)=\CO(n^3)$, by Theorem \ref{thm:main1}(b) (cf. the proof of (a) above). Similarly as above, Dynkin graph recognition does not affect the general complexity $\CO(n^3)$.\hfill$\Box$

 \section{Experiments and conclusions}

 In the experiments below we use our implementations of Algorithms \ref{alg:pdti} and \ref{alg:pdtri} in the computer algebra system Maple \cite{CIG}.

 We perform the test on an example coming from \lq\lq nature'', i.e., from practical application of Coxeter spectral analysis in representation theory of algebras, see \cite{Minfl2016,Mgroth2016}. We consider the following family of dense integral matrices.
 Let $\Nak{n}=[a_{i,j}]\in\MM_n(\ZZ)$ be the upper-triangular unidiagonal matrix whose coefficients are defined as follows:
 $$a_{i,i+s}:=\left\{\begin{array}{rr}
 \!\!1,&\mbox{if}\  s\, {\rm mod}\, 2 = 0,\\
 \!\!-1,&\mbox{if}\  s\, {\rm mod}\, 2 = 1,\end{array}\right.$$
 for any $1\leq i\leq n$ and $0\leq s\leq n-i$ (and $a_{i,j}=0$ for $i>j$).
 \begin{ex}\label{ex:Nakayama}
 For $n=2,4,5$, the matrices $\Nak{n}$ and the corresponding bigraphs $\Delta(\Nak{n})$ look as follows:
 \def\lewo{\!\!}
 $$\hspace{-0.6cm}\begin{array}{ccc}
 \Nak{2}&\Nak{4}&\Nak{5}\\
  {\scriptsize\left[ \begin {array}{rr} 1&\lewo-1\\ 0&\lewo1\end {array}
 \right]}&
 {\scriptsize\left[ \begin {array}{rrrr} 1&\lewo-1&\lewo1&\lewo-1\\ 0&\lewo1&\lewo-1&\lewo1\\ 0&\lewo0&\lewo1&\lewo-1\\ 0&\lewo0&\lewo0&\lewo1\end {array}
 \right]} &{\scriptsize\left[ \begin {array}{rrrrr}
 1&\lewo-1&\lewo1&\lewo-1&\lewo1\\ 0&\lewo1&\lewo-1&\lewo1&\lewo-1\\
 0&\lewo0&\lewo1&\lewo-1&\lewo1\\ 0&\lewo0&\lewo0&\lewo1&\lewo-1
\\ 0&\lewo0&\lewo0&\lewo0&\lewo1\end {array} \right]}\\
{\scriptsize
\xymatrixrowsep{15pt}
\xymatrixcolsep{15pt}
{\xymatrix{
  1\ar@{-}[d]  \\
2
}}}&{\scriptsize
\xymatrixrowsep{10pt}
\xymatrixcolsep{10pt}
{\xymatrix{
 & 1\ar@{-}[rd]\ar@{.}[dd]\ar@{-}[ld] &  \\
4\ar@{.}[rr]\ar@{-}[rd]&   & 2\ar@{-}[ld] \\
 & 3 &
}}}&{\scriptsize
\xymatrixrowsep{6pt}
\xymatrixcolsep{4pt}
{\xymatrix{
 & &1\ar@{.}[lld]\ar@{-}[lddd]\ar@{.}[rddd]\ar@{-}[rrd]& & \\
5\ar@{-}[rrrr]\ar@{-}[rdd]\ar@{.}[rrrdd]& & & &2\ar@{-}[ldd]\\
 & & & & \\
 &4\ar@{.}[rrruu]\ar@{-}[rr]& &3& \\
}}}
 \end{array}$$
 \end{ex}
 Note that the underlying graph $\overline{\Delta^{(n)}}$ of $\Delta^{(n)}:=\Delta(\Nak{n})$, is a complete graph, for every $n\geq 1$; in particular, $\Delta^{(n)}$ is connected. Moreover, one shows that each $\Delta^{(n)}$ is the so-called {\em Nakayama bigraph} of type $(n,2)$ introduced in \cite{Mgroth2016}. These bigraphs encode the $K$-theory of Nakayama algebras, an interesting class of finite-dimensional algebras recently studied by many authors, see
 \cite[Section 5]{Mgroth2016} for details. It follows from general theory that each $\Nak{n}$ is positive definite of Dynkin type $\AA_n$ (cf. \cite[Lemma 5.3]{Mgroth2016}), i.e., our inflation algorithms, applied to $\Delta^{(n)}$, return the Dynkin graph $\AA_n$, for every $n\geq 1$, cf. Remark \ref{rem:implem}(e).

We test the matrix $\Nak{400}$ (i.e., $400\times 400$ matrix) for positive definiteness in Maple ver. 15 on PC computer with Intel Core i5-7500 CPU 3.4GHz. First we apply to $\sym{\Nak{400}}$ an efficient implementation of the standard positive definiteness test based on Sylvester criterion and Gaussian elimination, briefly outlined in Introduction (the procedure {\tt PosDefTestByGaussElim} in \cite{CIG}). It returns {\tt true} after
$$35.562\, {\rm sec}.$$
We note that the execution time of the standard Maple routine {\tt gausselim} (for Gaussian elimination) on $\sym{\Nak{400}}$ is similar: $35.125$ sec. The procedure {\tt PosDefTestByRootInflations}$(\Nak{400},\stg)$ (Algorithm \ref{alg:pdtri}), returns {\tt true} after
$$\approx 0.969\, {\rm sec}.\ \ (\mbox{and}\ 599\ \mbox{inflations}),$$
for every strategy $\stg\in\{0,1,2,3\}$, since for $\Nak{400}$, the Dynkin graph $\AA_{400}$ emerges already after the line {\tt 2} (cf. Remark \ref{rem:implem}(a)). Whereas Algorithm \ref{alg:pdti} works with the speed depending on the strategy as follows:
$$
\begin{array}{c}
{\tt PosDefTestByInflations}(\Nak{400},\stg)\smallskip\\
{\scriptsize\begin{array}{r|rrrrrr}
\stg& 0& 1& 2& 2& 3& 3\\\hline
\#{\rm inflations}&39800  & 398 & 774 & 757 &3059  & 3016 \\
{\rm time\  in\  sec.}&837.72  & 6.11 & 6.54 & 6.45 & 587.65 & 568.22 \\
\end{array}}
\end{array}
$$
We presented two different executions for $\stg=2$, and also for $\stg=3$, since recall that these two strategies yield nondeterministic executions, with numbers of inflations (and hence the execution times) depending on the current behavior of the random number generator, see Section 4. Note that  the strategy $\stg=0$ appears to not be very efficient for the matrix $\Nak{400}$; however, the strategies $\stg=1,2$ work very nicely, not much slower than Algorithm \ref{alg:pdtri} above. A bottle-neck of our implementation of the strategy $\stg=3$ is the random selection of an edge among all dotted edges (see Section 4). We are working on a more efficient solution. On the other hand, it is a challenging task to construct a different, more sophisticated strategy for choosing dotted edges in Algorithm \ref{alg:random}, than our four simple strategies $\stg$.

The results of the above test are perhaps not breath-taking, but observe that our algorithms  are competing with Gaussian elimination, one of the simplest and most efficient basic matrix algorithms. And our algorithms are conceptually quite complex and they potentially have  several technical bottle-necks, not as easy to implement efficiently. The reader is referred to the file {\tt PosDefTestExperiments.mws}  in \cite{CIG} for the details of this and other tests, including positive definite
matrices of other Dynkin types and non-positive definite
ones (the readers having no access to Maple are referred
to the PDF version of the worksheet, also available in \cite{CIG}).

\smallskip

Concluding, we have constructed algorithmic positive definiteness tests of the following properties:
\begin{itemize}
\item their optimistic case is a positive definite matrix on the input (in contrast to the standard tests, cf. Introduction and Remark \ref{rem:implem}(c));
\item they work noticeably faster than the standard Gaussian elimination test on certain classes of matrices;
\item they can have nice applications, e.g. for the sets of integral matrices for which it is  a priori known that many of them are positive definite;
\item provide an example of a Las Vegas randomized algorithm (i.e., the behavior  varies from execution to execution on the same data, but the number of performed operations is always bounded by the quantity depending on the size of the input only); such algorithms are very interesting from the point of view of theoretical computer science;
\item a side effect of the algorithms is the computation of the Dynkin type of the quadratic form $q_A$ of a positive definite matrix $A$, encoding several additional information on $q_A$, see Remark \ref{rem:implem}(e);
\item they show a non-trivial practical application of deep theoretical techniques originally developed for different purposes (see \cite{Minfl2016, Mgroth2016, SimSIAM2013, KosFI2012, Ovsienko, BarotPena99}), and having a potential for further algorithmic studies.
\end{itemize}

Recently in \cite{PAR} the authors presented  an algorithm of the complexity $\CO(n^3)$ to compute the Dynkin types for a  larger (than our \utis's) class of integral matrices, called {\em symmetrizable quasi-Cartan matrices}, having origins in Lie theory. This algorithm can also be used as a positive definiteness test, it also uses inflations, but it is essentially different than ours. 
One can then say that the \lq\lq natural'' lower bound for the pessimistic complexity of such algorithms is $\CO(n^3)$, the complexity of the standard Gaussian elimination. Surprisingly, recent
 results from \cite{MM} on some new properties of the root systems allow to consider an algorithm to compute the Dynkin type and to test the positive definiteness for quasi-Cartan matrices, which has the complexity  $\CO(n^2)$ (that is, the same as simply reading the input matrix!).
We present this algorithm as well as some new applications of inflation techniques in the subsequent paper \cite{MM2}.

\end{document}